\documentclass[11pt]{article}

\usepackage{relsize}
\usepackage{amssymb}
\usepackage{amsthm}
\usepackage{mathtools}

\usepackage{stackengine}

\usepackage{afterpage}
\usepackage{makebox}

\DeclareMathSymbol{\alphld}{\mathalpha}{AMSa}{"45}
\DeclareMathSymbol{\Alpha}{\mathalpha}{operators}{"41}

\usepackage{setspace}
\usepackage[usenames,dvipsnames]{color}
\usepackage{xspace}
\usepackage{enumitem}

\newlength{\temp}

\newcommand{\medcup}{{\textstyle\bigcup}}
\newcommand{\medcap}{{\textstyle\bigcap}}

\newcommand{\mysub}[1]{\hbox{\smaller[3]{$#1$}}}

\newcommand{\sub}[2]{#1_{\hbox{\smaller[3]{$#2$}}}}

\newcommand{\subM}[1]{#1_{\mysub{M}}}

\let\leq=\leqslant
\let\le=\leqslant

\let\ge=\geqslant

\definecolor{dark-gray}{gray}{0.35}

\def \no#1#2#3 {{\bf #1} (#3), #2.}
\def \eds#1#2 {#1, #2.}

\definecolor{grey}{rgb}{0.9, 0.9, 0.9}
\newlength{\savefboxrule}
\newtheorem{mydef}{Definition} 
\newtheorem{mylemma}{Lemma} 
\newtheorem{myproposition}{Proposition} 
\newtheorem{mytheorem}{Theorem} 
\newtheorem{mycorollary}{Corollary} 
\newtheorem{myremarks}{Remarks} 
\newtheorem{myremark}[myremarks]{Remark} 
\newtheorem{myexample}{Example} 

\newcommand{\Vars}[1]{\mathrm{Vars}(#1)}

\newcommand{\COMMENT}[1]{}

\newcommand{\card}[1]{{\left|{#1}\right|}}

\newcommand{\defAs}{\coloneqq}

\newcommand{\Nats}{\mathbb{N}}

\newcommand{\Pow}{{\mathrm{pow}}}
\newcommand{\pow}[1]{\Pow({#1})}

\newcommand{\powAst}{\Pow^{\ast}}

\newcommand{\powast}[1]{\powAst({#1})}

\newcommand{\Nodes}{\mathcal{N}}

\newcommand{\st}{\,\texttt{|}\:}

\newcommand{\rk}{\hbox{\sf rk}\;}
\newcommand{\dom}{\hbox{\sf dom}}
\newcommand{\Places}{\mathcal{P}}
\newcommand{\TARGETS}{\mathcal{T}}
\newcommand{\Targets}[1]{\TARGETS({#1})}

\newcommand{\myphi}{\Phi}

\newcommand{\boldcalV}{\mbox{\boldmath$\mathcal{V}$}}

\newcommand{\HF}{\textnormal{\textsf{HF}}\xspace}
\newcommand{\MLS}{\textnormal{\textsf{MLS}}\xspace}

\newcommand{\MLSC}{\MLS\raisebox{.8pt}{\hspace{-1pt}$\times$}\xspace}
\newcommand{\MLSuC}{\MLS\raisebox{.9pt}{$\otimes$}\xspace}

\newcommand{\BST}{\textnormal{\textsf{BST}}}

\newcommand{\BSTuC}{\textnormal{\textsf{BST}}\xspace\raisebox{.9pt}{$\otimes$}\xspace}

\newcounter{instr}

\newcounter{instrb}

\newcommand{\commentout}[1]{}

\begin{document}


\title{Decidability of Extensions of Boolean Set Theory and the Limits of Explicit Graph Certificates}
\author{
Pietro Ursino \\
\emph{Dipartimento di Matematica e Informatica, Universit\`a di Catania}\\
\emph{Viale Andrea Doria 6, I-95125 Catania, Italy.}  \\
\mbox{E-mail:} \texttt{
pietro.ursino@unict.it}
}

\maketitle
\begin{abstract}
\noindent
We study the satisfiability problem for extensions of Boolean Set Theory
($\BST$) by set operators admitting suitable constructors, with particular
attention to the graph-based decision method using accessible fulfilling
$C$-graphs. The decidability construction for the class of $H$-operators is
presented in the preceding sections. We reconsider the proposed
$NP$-completeness argument based on reducing a fulfilling graph to
polynomially many places. We prove that this reduction cannot work in
general: already for the unordered Cartesian product $\otimes$, there is
a family of satisfiable $\BST^\otimes$-conjunctions containing $O(m^2)$
literals such that \emph{every} fulfilling $\otimes$-graph has at least
$2^m-1$ places. The lower bound is derived directly from the fulfillment
conditions and does not depend on the size of a particular model.
Consequently, explicit polynomial-size graph certificates cannot establish
an $NP$ upper bound for the whole class by this method. This does not
rule out polynomially verifiable succinct graph representations or other
$NP$ certificates; the general $NP$ membership question is not settled here.
\end{abstract}
\section*{Introduction}

$\BST$ is the quantifier-free fragment of set theory consisting of the propositional closure of atoms of the following types:
\[
x=y \cup z \/,  \quad x=y \setminus z\/,  \quad  x\neq y\/
\]
where $x,y,z,x_1\cdots x_n$ stand for set variables.

The first aim of the present paper consists in showing that
the technique developed in \cite{CU21} can be strengthen in order to reach the decidability for languages extending $\BST$ by allowing the use of literals of the type $\mathcal{C}(x_1,\cdots ,x_n)=x$, where $\mathcal{C}$ belongs to a special class of set operators, called $H$-operators (Definition \ref{Hoper}).
Roughly speaking, for a given set assignment to the variables $\mathfrak{I}(x_i)$  an $H$-operator builds the set $\mathcal{C}(x_1,\cdots ,x_n)$ piecewise using the venn regions of the input variables.
All the theories, obtained extending $\BST$ by allowing the use of an $H$-operator $C$, are decidable (Theorem \ref{Decidability}).

The proof of the above result (using Lemmas \ref{C-GraphFulfillsMyphi} and \ref{rev-otimesGraphFulfillsMyphi}) is a slight modifications of that developed by Cantone and Ursino in \cite{CU21}, therefore we put it in the Appendix.

Our second aim is to examine the complexity limitation of the explicit
graph-witness approach. Section~\ref{sec:graph-obstruction} provides a
satisfiable family of $\BST^\otimes$-conjunctions for which every fulfilling
accessible $\otimes$-graph must have exponentially many places. Thus the
pruning argument formerly used to infer $NP$ membership from small
explicit graphs does not establish the claimed upper bound; no
$NP$-completeness theorem for the entire class is asserted here.

$\MLSC$ is the acronym for MultiLevel Syllogistic ($\MLS$) extended with the Cartesian product operator $\times$.
We recall that $\MLS$ is the quantifier-free fragment of set theory involving the Boolean set operators $\cup$, $\cap$, and $\setminus$, and the equality and membership predicates.\\
Analogously, $\MLSuC$ denotes the extension of $\MLS$ with the  unordered Cartesian product operator $\otimes$, namely the set operator defined by $s \otimes t \coloneqq \big\{ \{u,v\} \st u \in s \wedge v \in t \big\}$.

$\BSTuC$ and \BST$\times$ can be obtained, respectively, from $\MLSuC$ and $\MLSC$ simply by dropping membership predicates.

Historically, the applications we deal with are related to languages which are $\MLS$-extensions, but, in many cases, the $\in$-literals have no role, in the sense that the application runs in the same way without them. Therefore we can limit ourselves to detect the problems to their \BST-versions.

The well-celebrated Hilbert's Tenth problem (HTP, for short)\cite{Hilbert-02}, posed by David Hilbert at the beginning of last century, asks for a uniform procedure that can determine in a finite number of steps whether any given Diophantine polynomial equation with integral coefficients is solvable in the integers.

In 1970, it was shown that no algorithmic procedure exists for the solvability problem of generic polynomial Diophantine equations, as a result of the combined efforts of M.\ Davis, H.\ Putnam, J.\ Robinson, and Y.\ Matiyasevich (DPRM theorem, see \cite{Rob,DPR61,Mat70}).

In the early eighties, D.Cantone asked whether the decision problems for the theories $\MLSC$ and $\MLSuC$ can be reducible to HTP.
The decision problem for $\MLSC$ (resp., $\MLSuC$) can be somehow regarded as the set-theoretic counterparts of HTP, where the union of disjoint sets and the Cartesian product (resp., unordered Cartesian product of disjoint sets) play the roles of integer addition and multiplication, respectively.
Indeed, $|s \cup t| = |s| + |t|$, for any disjoint sets $s$ and $t$, and $|s \times t| = |s| \cdot |t|$, for any sets $s$ and $t$ (whereas $|s \otimes t| = |s| \cdot |t|$, for any \emph{disjoint} sets $s$ and $t$).

When $\MLSuC$ and $\MLSC$ are extended with the two-place predicate $|\cdot| \leq |\cdot|$ for cardinality comparison, where $|s| \leq |t|$ holds if and only if the cardinality of $s$ does not exceed that of $t$, the satisfiability problem for such extensions  become undecidable, since HTP would be reducible to each of them, as proved in \cite{CCP90}.
In the above cited reduction to HTP, membership operator plays no role \cite{CU18}, then, by extending $\BSTuC$ or \BST$\times$with the two-place predicate $|\cdot| \leq |\cdot|$, you get again a problem reducible to HTP.

Therefore, the real set-theoretic counterpart of HTP is actually $\BSTuC$ or $\BST\times$.

The presence of $\otimes$ and $\times$ among the operators under
consideration makes the distinction between decidability and sharp
complexity bounds essential. In particular, a decision procedure based
on finite graph witnesses need not yield polynomial-size witnesses.
Section~\ref{sec:graph-obstruction} establishes an exponential lower
bound for the explicit $\otimes$-graph representation.

\medskip

Notice that there are satisfiable $\BSTuC$ and $\BST\times$ formulas that admit only infinite models (for example $\myphi :=x\neq\emptyset\wedge x\supseteq x\otimes x$).
Nevertheless, we prove in \cite{CU14} that even theories which force a model to be infinite can be proven to be decidable by using the small witness-model property, which is a way to finitely represent the infinity.
Hence, the decidability of $\BSTuC$ or $\BST\times$ cannot be proven through a small model property. In \cite{CU21} it is proven that $\BSTuC$ is decidable by way of an algorithmic representation.
Our general results prove, by-product, that $\BST\times$ behaves in the same way.

In view of the above observation, to be more precise, the set-theoretic counterpart to HTP is the \emph{finite} satisfiability problem for $\BSTuC$ or $\BST\times$, namely the problem of establishing algorithmically for any $\BSTuC$($\BST\times$)-formula whether it admits a finite model or not.

In \cite{CU21} we use a special class of algorithms to decide this problem.
In the same fashion the proof of Theorem \ref{Decidability} runs for finite satisfiability of \BST$\mathcal{C}$.

No $NP$ upper bound for $\BST^\otimes$, $\BST\times$, or
$\BST P$ follows from the graph-pruning argument reconsidered below.
The existing decidability conclusions should therefore be kept separate
from complexity claims. For context, Cantone and Maugeri
\cite{CM23Tensor} establish $NP$-completeness for a restricted fragment
of $\BST^\otimes$ with only logarithmically many distinct left-hand-side
variables of product literals; this does not resolve the unrestricted case.

As shown in \cite{Schw78,CCS90}, the finite satisfiability property for the extension of \MLS with cardinality comparison holds if and only if the cardinality of $s$ does not exceed that of $t$, can be reduced to purely existential Presburger arithmetic, which is known to be NP-complete (see \cite{Sca}).
The contrast between these arithmetic reductions and decidable set
fragments remains of interest, but does not by itself yield an
$NP$-completeness statement for set languages with product.

\section{Extensions of $\BST$ with $\mathcal{C}$ operator}\label{Assignment}

By "Sets" we mean the von Neumann universe $\boldcalV$ of all well-founded sets.

We recall that $\boldcalV$ is a cumulative hierarchy constructed in stages by transfinite recursion over the class $\mathit{On}$ of all ordinals. Specifically, $\boldcalV \defAs \medcup_{\alpha \in \mathit{On}} \boldcalV_{\alpha}$ where, recursively, $\boldcalV_{\alpha} \defAs \medcup_{\beta<\alpha} \pow{\boldcalV_{\beta}}$, for every $\alpha \in \mathit{On}$, with $\pow{\cdot}$ denoting the powerset operator. Based on that construction, we can readily define the \textsc{rank} of any well-founded set $s \in \boldcalV$, denoted $\rk{s}$, as the least ordinal $\alpha$ such that $s \subseteq \boldcalV_{\alpha}$. The set of the sets of finite rank, hence belonging to $\boldcalV_{\alpha}$ for some finite ordinal $\alpha$, forms the set \HF of the \textsc{hereditarily finite sets}. Thus, $\HF = \boldcalV_{\omega}$, where $\omega$ is the first limit ordinal, namely the smallest non-null ordinal with no immediate predecessor.

$\BST\mathcal{C} $ is the quantifier-free fragment of set theory consisting of the propositional closure of atoms of the following types:
\[
x=y \cup z \/,  \quad x=y \setminus z\/,  \quad x = \mathcal{C}(x_1\cdots x_n) \/,   \quad  x\neq y\/
\]
where $x,y,z,x_1\cdots x_n$ stand for set variables and $\mathcal{C}$ is an unordered (ordered) operator on sets, that is a map $\mathcal{C}(x_1\cdots x_n):\mbox{Sets}\otimes\cdots\otimes\mbox{Sets}\rightarrow\mbox{Sets}$ (if ordered, $\mathcal{C}:\mbox{Sets}\times\cdots\times\mbox{Sets}\rightarrow\mbox{Sets}$).

For any $\BST\mathcal{C} $-formula $\myphi$, we denote by $\Vars{\myphi}$ the collection of set variables occurring in it.

The semantics of $\BST\mathcal{C}$ follows exactly the definition of many other languages already treated (see for example \cite{CU18,CU21}). Other operators and relators of $\BST\mathcal{C} $ are interpreted according to their usual semantics as well as satisfiability by partitions and normalization of $\BST\mathcal{C}$-formulae \cite{CU21}.

You can find all the following definitions in previous articles such as \cite{CU14,CU18,CU21}.
For the reader's convenience we briefly resume them in the following section.

\subsection{Satisfiability by partitions}

A \textsc{partition} is a collection of pairwise disjoint non-null sets, called the \textsc{blocks} of the partition. The union $\bigcup \Sigma$ of a partition $\Sigma$ is its \textsc{domain}.

A \textsc{set assignment} $M$ is any map from a collection $V$ of set variables (called the \textsc{variable-domain of $M$} and denoted $\dom(M)$) into the von Neumann universe $\boldcalV$ of all well-founded sets.

Let $V$ be a finite collection of set variables and $\Sigma$ a  partition. Also, let $\mathfrak{I} \colon V \rightarrow \pow{\Sigma}$ be any map. In a very natural way, the map $\mathfrak{I}$ induces a set assignment $M_{\mysub{\mathfrak{I}}}$ over $V$ defined by:
\[
\textstyle
M_{\mysub{\mathfrak{I}}} v \defAs \bigcup \mathfrak{I}(v)\/, \qquad \text{for $v \in V$\/.}
\]
We refer to the map $\mathfrak{I}$ (or to the pair $(\Sigma, \mathfrak{I})$, when we want to emphasize the partition $\Sigma$) as a \textsc{partition assignment}.

Whenever $\Sigma^*\subseteq\Sigma$, it arises automatically an assignment restricted to $\Sigma^*$ defined as follow $$\mathfrak{I}^*(x)=\mathfrak{I}(x)\cap\Sigma^*$$.

\begin{mydef}\label{def:satisfiability}\rm
Let $\Sigma$ be a partition and $\mathfrak{I} \colon V \rightarrow \pow{\Sigma}$ be a partition assignment over a finite collection $V$ of set variables. Given a $\BST \mathcal{C}$-formula $\myphi$ such that $\Vars{\myphi} \subseteq V$, we say that $\mathfrak{I}$ \textsc{satisfies $\myphi$}, and write $\mathfrak{I} \models \myphi$, when the set assignment $M_{\mysub{\mathfrak{I}}}$ induced by $\mathfrak{I}$ satisfies $\myphi$ (equivalently, one may say that $\Sigma$ \textsc{satisfies $\myphi$ via the map $\mathfrak{I}$}, and write $\Sigma/\mathfrak{I} \models \myphi$, if we want to emphasize the partition $\Sigma$). We say that $\Sigma$ \textsc{satisfies} $\myphi$, and write $\Sigma \models \myphi$, if $\Sigma$ satisfies $\myphi$ via some map $\mathfrak{I} \colon  V \rightarrow \pow{\Sigma}$.
\end{mydef}

The following result can be proved immediately.
\begin{mylemma}\label{wasB}
If a \BST$\mathcal{C}$-formula is satisfied by a partition $\Sigma$, then it is satisfied by any partition $\overline \Sigma$ that includes $\Sigma$ as a subset, namely such that $\Sigma \subseteq \overline \Sigma$.
\end{mylemma}

Plainly, a \BST$\mathcal{C}$-formula $\myphi$ satisfied by some partition is satisfied by a set assignment. Indeed, if $\Sigma \models \myphi$, then $\Sigma/\mathfrak{I} \models \myphi$ for some map $\mathfrak{I} \colon  V \rightarrow \pow{\Sigma}$, and therefore $M_{\mysub{\mathfrak{I}}} \models \myphi$.
The converse holds too. In fact, let us assume that $M \models \myphi$, for some set\index{set} assignment $M$ over the collection $V = \Vars{\myphi}$ of the set variables occurring in $\myphi$, and let $\subM{\Sigma}$ be the \textsc{Venn partition} induced by $M$, namely
\[
\sub{\Sigma}{M} \defAs \Big\{ \medcap MV' \setminus \medcup M(V \setminus V') \st \emptyset \neq V' \subseteq V \Big\} \setminus \big\{ \,\emptyset\,\big\}.\footnotemark
\]
\footnotetext{Hence, we have:
\begin{enumerate}[label=-]
\item $(\forall \sigma \in \sub{\Sigma}{M})(\forall v \in V)(\sigma \cap Mv = \emptyset \vee \sigma \subseteq Mv)$,

\item $(\forall \sigma,\sigma' \in \sub{\Sigma}{M}) \big( (\forall v \in V) (\sigma \subseteq Mv \leftrightarrow \sigma' \subseteq Mv) \leftrightarrow \sigma = \sigma' \big)$, and

\item $\medcup \Sigma = \medcup MV$.
\end{enumerate}
}
Let $\subM{\mathfrak{I}} \colon V \rightarrow \pow{\subM{\Sigma}}$ be the map defined by
\[
\subM{\mathfrak{I}}(v) \defAs \{ \sigma \in \subM{\Sigma} \st \sigma \subseteq Mv\}\/, \qquad \text{for $v \in V$.}
\]
It is an easy matter to check that the set assignment induced by $\subM{\mathfrak{I}}$ is just $M$. Thus $\subM{\Sigma}/\subM{\mathfrak{I}} \models \myphi$, and therefore $\subM{\Sigma} \models \myphi$, proves that $\myphi$ is satisfied by some partition, in fact by the Venn partition induced by $M$, whose size is at most $2^{|V|}-1$.

Thus, the notion of satisfiability by set assignments and that of satisfiability by partitions coincide.

As a by-product of Lemma~\ref{wasB} and the above considerations, we also have:
\begin{mylemma}\label{wasA}
Every \BST$\mathcal{C}$-formula $\myphi$ with $n$ distinct variables is satisfiable if and only if it is satisfied by some partition with $2^{n}-1$ blocks.
\end{mylemma}

Regarding the literals of the type $\{\cup,\setminus\}$ the assignment models a literal whenever, applied to the partition, verifies the literal \cite{CU18,CU21}.
\begin{mylemma}\label{partitionAssignmentBoolean}
Let $\Sigma$ be a partition and let $\mathfrak{I} \colon V \rightarrow \pow{\Sigma}$ be a partition assignment over a (finite) set of variables $V$. Then, for all $x,y,z \in V$ and $\star \in \{\cup,\setminus\}$, we have:
\begin{enumerate}[label=(\alph*)]
\item\label{partitionAssignmentBooleanA} $\mathfrak{I} \models x=y \star z \quad \Longleftrightarrow \quad \mathfrak{I}(x) = \mathfrak{I}(y) \star \mathfrak{I}(z)$,

\item\label{partitionAssignmentBooleanB} $\mathfrak{I} \models x\neq y\quad \Longleftrightarrow \quad \mathfrak{I}(x) \neq \mathfrak{I}(y)$.
\end{enumerate}
\end{mylemma}
\subsection{Normalization of \BST$\mathcal{C}$-formulae}

By applying disjoint normal form and the simplification rules illustrated in \cite{CU18}, the satisfiability problem for \BST$\mathcal{C}$ can be reduced to the satisfiability problem for \textsc{normalized conjunctions} of \BST$\mathcal{C}$, namely
conjunctions of \BSTuC-literals of the following restricted types:%
\begin{gather}\label{formula}
  x=y \cup z \/,\ \  x=y \setminus z\/,\ \    x = \mathcal{C}(x_1\dots x_n)\/,\ \  x \neq y\/,
  \end{gather}
where $x,y,z$ stand for set variables. Indeed, it is enough to observe that:
\begin{itemize}
\renewcommand\labelitemi{--}

\item $x \subseteq y$ is equivalent to $x = x \cap y$;

\item the terms $y \cap z$ and $y \setminus (y \setminus z)$ are equivalent, so  an atom of the form $x = y \cap z$ is equisatisfiable with the conjunction $x = y \setminus y' \wedge y' = y \setminus z$, where $y'$ stands for any fresh set variable;

\item each negative literal of the form $x \neq y \star z$ (with $\star \in \{\cup, \cap,\setminus,\otimes\}$) is equisatisfiable with the conjunction $x' = y \star z \wedge x' \neq x$, where $x'$ stands for any fresh set variable.
\end{itemize}

We shall refer to literals of the form $x = \mathcal{C}(x_1\dots x_n)$ as \textsc{$\mathcal{C}$-literals}, and similarly for the other types of literals in \eqref{formula}.

\section{$C$-graphs}
In this section we show in details how to build $\mathcal{C}(x_1,\dots ,x_n)$ from the Venn regions  of $x_1,\dots ,x_n$. First we have to define a set-constructor on a partition.

An unordered constructor on a partition $\Sigma$ is a map $C:pow(\Sigma)\rightarrow \ Sets$ (for example the unordered cartesian product $\otimes$).
An ordered constructor on a partition $\Sigma$ is a map $C:Seq(\Sigma)\rightarrow \ Sets$, where $Seq(\Sigma)$ are finite ordered sequences $<\sigma_1\cdots\sigma_k>$ of blocks $\sigma_i\in\Sigma$ (for example the ordered cartesian product $\times$).
In a natural way all the concepts and results developed in Section \ref{Assignment} can be generalized to $Seq(\Sigma)$.
\begin{mydef}\label{disjoint}
  Let $\Sigma$ be a partition. A set constructor $C:\pow\Sigma\rightarrow\mbox{ Sets }$ with input an ordered or unordered subset of $\Sigma$ is {\bf disjoint} if for all $X,Y\subset\Sigma$ $X\neq Y$, $C(X)\cap C(Y)=\emptyset$ (if ordered $X,Y\in Seq(\Sigma$)).
\end{mydef}

\noindent
From now on our constructors are supposed to be disjoint.
At this point we have to set which subsets of $\Sigma$ are needed in order to construct $\mathcal{C}(x_1,\dots ,x_n)$. $Q$ is meant to be the function that selects such subsets.

Let $\mathcal{C}$ be an unordered (ordered) operator of -arity $n$. If for any
set of variables $V$, any partition $\Sigma$ and any partition assignment $\mathfrak{I}:V\rightarrow pow(\Sigma)$ ($\mathfrak{I}:V\rightarrow Seq(\Sigma)$) there exists a function $Q^{(\mathcal{C},\Sigma ,\mathfrak{I})}:\underbrace{V\otimes\cdots\otimes V}_n\rightarrow pow(pow(\Sigma))$ ($Q^{((\mathcal{C},\Sigma ,\mathfrak{I})}:\underbrace{V\times\cdots\times V}_n\rightarrow pow(Seq(\Sigma))$) and a unordered (ordered) constructor $C^{(\mathcal{C},\Sigma ,\mathfrak{I})}$ such that
\vspace{0.3cm}
\begin{description}
  \item[(a)] $\bigcup_{N\in Q^{(\mathcal{C},\Sigma ,\mathfrak{I})}(x_1,\dots,x_n)} C^{(\mathcal{C},\Sigma ,\mathfrak{I})}(N)=\mathcal{C}(\medcup\mathfrak{I}(x_1)\cdots\medcup\mathfrak{I}(x_n)),\ \bigcup\{N\st N\in Q^{((\mathcal{C},\Sigma ,\mathfrak{I})}(x_1,\dots,x_n)\}=\mathfrak{I}(x),$
      {\color{red} forse \'e target di N ma molto probabilmente non \'e necessaria}
  \item[(b)] For any subpartition $\Sigma^*$ of $\Sigma$ and any $x_1,\dots,x_n\in V$, $Q^{(\mathcal{C},\Sigma^* ,\mathfrak{I}\mid_{\Sigma^*})}(x_1,\dots,x_n)=Q^{(\mathcal{C},\Sigma ,\mathfrak{I})}(x_1,\dots,x_n)\}\cap \Pow\Sigma^*$ and $C^{(\mathcal{C},\Sigma^* ,\mathfrak{I}\mid_{\Sigma^*})}=C^{(\mathcal{C},\Sigma ,\mathfrak{I})}\mid_{\Sigma^*}$.
\end{description}

By $Q^{(\mathcal{C},\Sigma ,\mathfrak{I})}(x_1,\dots,x_n)$ we are briefly indicating the image of $Q^{(\mathcal{C},\Sigma ,\mathfrak{I})}(x_1,\dots,x_n)$.

\vspace{0.3cm}

\noindent
then $\mathcal{C}$ is an unordered (ordered) operator \textsc{with constructors} $(C,Q)$.

\begin{myexample}\label{constructors}

We now give some examples of operators with constructors.

\begin{myproposition}
 $\mathcal{C}=\Pow,\times,\otimes$ are operators with constructors.
\end{myproposition}
 \begin{proof}

In all cases, we assume that a partition $\Sigma$ and a partition assignment $\mathfrak{I}:V\rightarrow pow(\Sigma)$ (referred to the involved variables), are given.

For each operator $\Pow,\times,\otimes$ we are going to find the following couples, respectively

   $(C^{(\Pow ,\Sigma ,\mathfrak{I})},Q^{(\Pow ,\Sigma ,\mathfrak{I})})$,
   $(C^{(\times ,\Sigma ,\mathfrak{I})},Q^{(\times ,\Sigma ,\mathfrak{I})})$,
   $(C^{(\otimes ,\Sigma ,\mathfrak{I})},Q^{(\otimes ,\Sigma ,\mathfrak{I})})$.

   In order to define $C^{(\Pow ,\Sigma ,\mathfrak{I})}$ we introduce a modification of powerset operator, the operator $\powAst$, that has, as argument, a collection of sets $S$, and it is defined as follow

   \begin{align*}
\powast{S} &\defAs \big\{t\mid t \subseteq \medcup S \land t \cap s \neq \emptyset, \text{ for every } s \in S \big\}
\end{align*}
then
\begin{align*}
C^{(\Pow ,\Sigma ,\mathfrak{I})}(S)=\powast{S}
\end{align*}

\begin{align*}
C^{(\times ,\Sigma ,\mathfrak{I})}(\sigma_1,\sigma_2)=\sigma_1\times\sigma_2 &\defAs \big\{<x,y> \mid x \in \sigma_1,\ y \in \sigma_2, <x,y>\text{ ordered couple } \big\}
\end{align*}

\begin{align*}
C^{(\otimes ,\Sigma ,\mathfrak{I})}(\sigma_1,\sigma_2)=\sigma_1\otimes\sigma_2 &\defAs \big\{\{x,y\} \mid x \in \sigma_1,\ y \in \sigma_2 \big\}
\end{align*}

Regarding $Q$:

For a literal $pow(x)=y$; $Q^{(\Pow ,\Sigma ,\mathfrak{I})}(x)= \{\{q_1\cdots q_n\}\mid\{q_1\cdots q_n\}\subseteq\mathfrak{I}(x) \}$ (i.e. $\Pow\mathfrak{I}(x)$).

For a literal $x\times y=z$; $Q^{(\times ,\Sigma ,\mathfrak{I})}(x,y)= \{<q_1,q_2>\mid q_1\in\mathfrak{I}(x)\ q_2\in\mathfrak{I}(y)\}$(i.e. $\mathfrak{I}(x)\times\mathfrak{I}(y)$).

For a literal $x\otimes y=z$; $Q^{(\otimes ,\Sigma ,\mathfrak{I})}(x,y)= \{\{q_1,q_2\}\mid q_1\in\mathfrak{I}(x)\ q_2\in\mathfrak{I}(y)\}$(i.e. $\mathfrak{I}(x)\otimes\mathfrak{I}(y)$).

The verification of condition (a) is a straight calculation of the following equalities.

$$\bigcup_{N\in Q^{(\Pow ,\Sigma ,\mathfrak{I})}(x)}C^{(\Pow ,\Sigma ,\mathfrak{I})} {N}=\Pow M_{\mysub{\mathfrak{I}}}x$$
$$\bigcup_{N\in Q^{(\times ,\Sigma ,\mathfrak{I})}(x,y)}C^{(\times ,\Sigma ,\mathfrak{I})} {N}=M_{\mysub{\mathfrak{I}}}x\times M_{\mysub{\mathfrak{I}}}y$$
$$\bigcup_{N\in Q^{(\otimes ,\Sigma ,\mathfrak{I})}(x,y)}C^{(\otimes ,\Sigma ,\mathfrak{I})} {N}=M_{\mysub{\mathfrak{I}}}x\otimes M_{\mysub{\mathfrak{I}}}y$$

Regarding the condition (b) of $\Pow$,

$$Q^{(\Pow ,\Sigma^* ,\mathfrak{I}^*)}(x)=\Pow\mathfrak{I}(x)\cap\\Pow\Sigma^*=Q^{(\Pow ,\Sigma ,\mathfrak{I})}(x)\cap\\Pow\Sigma^*$$.

Regarding the condition (b) of $\times$,  $$Q^{(\times ,\Sigma^* ,\mathfrak{I}^*)}(x,y)=\mathfrak{I}^*(x,y)\times \mathfrak{I}^*(y)=\mathfrak{I}(x)\times \mathfrak{I}(y)\cap\Sigma^*\times\Sigma^*=\mathfrak{I}(x)\times \mathfrak{I}(y)\cap Seq\Sigma^*=Q^{(\times ,\Sigma ,\mathfrak{I})}(x,y)\cap Seq\Sigma^*$$.

Regarding the condition (b) of $\otimes$, $$Q^{(\otimes ,\Sigma^* ,\mathfrak{I}^*)}(x,y)=\mathfrak{I}^*(x)\otimes \mathfrak{I}^*(y)=\mathfrak{I}(x)\otimes \mathfrak{I}(y)\cap\Sigma^*\otimes\Sigma^*=\mathfrak{I}(x)\otimes \mathfrak{I}(y)\cap\Pow\Sigma^*=Q^{(\otimes ,\Sigma ,\mathfrak{I})}(x,y)\cap \Pow\Sigma^*$$.

 \end{proof}
\end{myexample}

\noindent
We call informally $C$-products the sets built using the constructor $C$.

Roughly speaking, the main aim of the following definition is to detect which blocks receives $C$-products, and we will call them $\sub{\Sigma}{C}$,
and which sets of blocks creates $C$-products, and we will call them $\sub{\Pi}{C}$. For example, if $C$ is $\otimes$, inside $\sub{\Sigma}{\otimes}$ there are blocks composed only of pairs $\{s,t\}$ with $s\in\sigma$ and $t\in\sigma'$ where $\sigma ,\sigma'$ are blocks of the given partition $\Sigma$. Instead, $\sub{\Pi}{\otimes}$ is composed exactly of that couples of blocks, $\{\sigma ,\sigma'\}$ that creates the pairs.

\begin{mydef}\label{C-blocks-upblocks}\rm
A subset $\Sigma^{*}$ of a partition $\Sigma$ is a \textsc{$C$-subpartition} of $\Sigma$  if $\medcup \Sigma^{*} = \medcup C[\mathcal{B}]$, for some $\mathcal{B} \subseteq Q^{(C,\Sigma)}$. By abuse of notation,f or a generic collection of variables and a set assignment on $\Sigma$, $\mathfrak{I}$, we denote by $Q^{(\Sigma)}$ the image of $Q^{(\mathcal{C},\Sigma ,\mathfrak{I})}$

We denote by $\sub{\Sigma}{C}$ the $\subseteq$-maximal $C$-subpartition of $\Sigma$ and we refer to its elements as the \textsc{$C$-blocks} of $\Sigma$.
We also denote by $\sub{\Pi}{C}$ the subset of $pow(\Sigma)$ such that $\medcup \sub{\Sigma}{C} = \medcup C[\sub{\Pi}{C}]$ and we refer to its elements as \textsc{$C$-upblocks} (`upblocks' for \emph{subsets of blocks which generates through the constructor $C$ the partition $\sub{\Sigma}{C}$}).
\end{mydef}

\begin{mydef}[$C$-graphs] \rm
A \textsc{$C$-graph} $\mathcal{G}$ is a directed bipartite graph whose set of vertices comprises two disjoint parts: a set of \textsc{places} $\Places$ and a set of \textsc{nodes} $\Nodes$, where $\Nodes = \pow\Places$ (if $C$ is ordered $\Nodes = Seq(\Places$)). The edges issuing from each place $q$ are exactly all pairs $\langle q,B \rangle$ such that $q\in B \in \sub{\Nodes}{}$\/: these are called \textsc{membership edges}.
The remaining edges of $\sub{\mathcal{G}}{}$, called \textsc{distribution} or \textsc{saturation edges}, go from nodes to places. When there is an edge $\langle B,q \rangle$ from a node $B$ to a place $q$, we say that $q$ is a \textsc{target} of $B$.
The map $\sub{\TARGETS}{}$ over $\sub{\Nodes}{}$ defined by
\[
\sub{\TARGETS}{}(B) \defAs \{q \in \Places \st q \text{ is a target of } B\}, \quad \text{for } B \in \sub{\Nodes}{},
\]
is the \textsc{target map} of $\sub{\mathcal{G}}{}$.

As soon as we will show the semantic counterpart of $\mathcal{G}$, it will be clear that the targets are the candidate to verifies $\TARGETS [\sub{\Nodes}{}]\subseteq \sub{\Sigma}{C}$.
The \textsc{size} of $\sub{\mathcal{G}}{}$ is defined as the cardinality of its set of places $\sub{\Places}{}$. Plainly, a $C$-graph $\sub{\mathcal{G}}{}$ is fully characterized by its target map $\sub{\TARGETS}{}$, since the sets of nodes and of places of $\sub{\mathcal{G}}{}$ are expressible as $\dom(\TARGETS)$ and $\medcup \dom(\TARGETS)$, respectively. When convenient, we shall explicitly write $\sub{\mathcal{G}}{} = (\Places, \Nodes, \TARGETS)$ for a $C$-graph with set of places $\Places$, set of nodes $\Nodes$, and target map $\TARGETS$ (in order to simplify the graph, if $N\in\Nodes$ has not outgoing arrows we drop $N$).
\end{mydef}

\smallskip
We now give a semantic version of the definition of $C$-graphs. For this purpose,
we show how to construct the $C$-graph $\sub{\mathcal{G}}{\Sigma}$ \textsc{induced by a given partition} $\Sigma$.

To begin with, we select a set of \textsc{places} $\sub{\Places}{\Sigma}$ of the same cardinality of $\Sigma$ such that $\sub{\Places}{\Sigma}$ and $pow(\sub{\Places}{\Sigma})$ are disjoint, and define the vertex set of $\sub{\mathcal{G}}{\Sigma}$ as the union $\sub{\Places}{\Sigma} \cup (pow(\sub{\Places}{\Sigma}))$.
The members of $pow(\sub{\Places}{\Sigma})$  will form the set of \textsc{nodes} $\sub{\Nodes}{\Sigma}$ of $\sub{\mathcal{G}}{\Sigma}$.
Places in $\sub{\Places}{\Sigma}$ are intended to be an abstract representation of the blocks of $\Sigma$ via a bijection $q \mapsto q^{(\bullet)}$ from $\sub{\Places}{\Sigma}$ onto $\Sigma$.
Likewise, nodes in $\sub{\Nodes}{\Sigma}$ are intended to represent the application of the constructor $C$ to the blocks represented by their places.

The disjoint sets $\sub{\Places}{\Sigma}$ and $\sub{\Nodes}{\Sigma}$ will form the parts of the bipartite graph $\sub{\mathcal{G}}{\Sigma}$ we are after.
The bijection $(\bullet)$ can be naturally extended to nodes $B$ of $\sub{\mathcal{G}}{\Sigma}$ by putting $B^{(\bullet)} \defAs \{q^{(\bullet)} \st q \in B\}$.

Having defined the vertex set of $\sub{\mathcal{G}}{\Sigma}$, next we describe its edge set.
The edges issuing from each place $q$ are exactly all pairs $\langle q,B \rangle$ such that $q\in B \in \sub{\Nodes}{\Sigma}$ (membership edges of $\sub{\mathcal{G}}{\Sigma}$).
The remaining edges of $\sub{\mathcal{G}}{\Sigma}$ go from nodes to places (distribution or saturation edges of $\sub{\mathcal{G}}{\Sigma}$).
Only places $q$ corresponding to $C$-blocks $q^{(\bullet)}$ of $\Sigma$ (hence called $C$-places) can have incoming edges. Likewise, only nodes $B$ such that $B^{(\bullet)} \in \sub{\Pi}{\mathcal{C}}$ (the candidate to be the set $Q^{(\Sigma)}$ of the constructible operator) can have outgoing edges. Such nodes will be called $Q^{(\Sigma)}$-nodes.
Specifically, for a $Q^{(\Sigma)}$-node $B$ and a $C$-place $q$ of $\sub{\mathcal{G}}{\Sigma}$, there is an edge $\langle B,q \rangle$ exactly when
\[
q^{(\bullet)} \cap C(B^{(\bullet)}) \neq \emptyset,
\]
 namely when there is some ``flow'' of elements built by constructor $C$ applied to the node $B^{(\bullet)}$
 from $C(B^{(\bullet)})$ to $q^{(\bullet)}$ (through the edge $\langle B,q \rangle$). This is the sense in which a $C$-graph can be considered a kind of flow graph in the realm of set theory.
Thus, the target map $\sub{\TARGETS}{\Sigma}$ of $\sub{\mathcal{G}}{\Sigma}$ is defined by
\[
\sub{\TARGETS}{\Sigma}(B) \defAs \{q \in \sub{\Places}{\Sigma,\! C} \st q^{(\bullet)} \cap \mathcal{C}(B^{(\bullet)}) \neq \emptyset\}, \quad \text{for } B \in \sub{\Nodes}{\Sigma,\! C},
\]
where $\sub{\Places}{\Sigma,\! C}$ and $\sub{\Nodes}{\Sigma,\! C}$ denote the collections of the $C$-places and of the $Q^{(\Sigma)}$-nodes of $\sub{\mathcal{G}}{\Sigma}$, respectively.

Notice that each $Q^{(\Sigma)}$-node $B$ of $\sub{\mathcal{G}}{\Sigma}$ has some target. Indeed, from $\medcup \sub{\Sigma}{C} = \medcup C[\sub{\Pi}{\mathcal{C}}]$ it follows that $\emptyset \neq C(B^{(\bullet)}) \subseteq \medcup \sub{\Sigma}{C}$, and therefore $\sub{\TARGETS}{\Sigma}(B) \neq \emptyset$.

\begin{myexample}\label{cgraph}
  In order to make more clear the above definitions and notations, we present a concrete example.

  Consider the following unordered constructor:

  $$2s\otimes(\sigma_1,\sigma_2)=\{\{\{x'\},\{y'\}\}\}\mid x'\in\sigma_1,y'\in\sigma_2\}.$$

  Consider $\mathcal{C}=2s\otimes$ ({\bf Two Singlelizer operator}) so defined $$2s\otimes(x,y)=\{\{\{x'\},\{y'\}\}\mid x'\in x,y'\in y\}.$$

  Let $\Sigma$ a partition, $V$ a set of variables and $\mathfrak{I}$ a set assignment to the variables. For a literal of the type $2s\otimes$, $Q^{(2s\otimes,\Sigma,\mathfrak{I})}(x,y)$ is equal to $\{\{\sigma_1,\sigma_2\}\mid \sigma_1\in\mathfrak{I}(x),\sigma_2\in\mathfrak{I}(y)\}\}$.

  It can be seen that $2s\otimes$ is an operator with constructors $(C^{(2s\otimes,\Sigma,\mathfrak{I})},Q^{(2s\otimes,\Sigma,\mathfrak{I})})$.

  In this case, for abuse of notation, we identify $C$ and $\mathcal{C}$ with $2s\otimes$ (in general the constructor could be different from the operator for example $\Pow$).

  Denote $\emptyset^1=\{\emptyset\},\emptyset^2=\{\{\emptyset\}\}\dots\emptyset^n=\underbrace{\{\dots\{}_n\emptyset\}\dots\}$

  Let $\sigma_1'^{(\bullet)}=\{\emptyset^1\},\sigma_1''^{(\bullet)}=\{\emptyset^4\},\sigma_2=^{(\bullet)}\{\emptyset^2\} $
  then

  $$2s\otimes(\sigma_1',\sigma_2)\cup 2s\otimes(\sigma_1'',\sigma_2)=\{\{\emptyset^2,\emptyset^3\},\{\emptyset^5,\emptyset^3\}\}$$.

   Consider the following set of variables $V:=\{x,y_1,y_2\}$, the partition $\Places_\Sigma=\{\sigma_1',\sigma_1'',\sigma_2,\sigma_5',\sigma_5''\}$ with
    $\sigma_5'^{(\bullet)}=\{\{\emptyset^2,\emptyset^3\}\}$ and
      $\sigma_5''^{(\bullet)}=\{\{\emptyset^5,\emptyset^3\}\}$,
      and the partition assignment $\mathfrak{I}:V\rightarrow pow(\Sigma)$ so defined $\mathfrak{I}(x)=\{\sigma_5',\sigma_5''\},\mathfrak{I}(y_1)=\{\sigma_1',\sigma_1''\},
      \mathfrak{I}(y_2)=\{\sigma_2\}$.
      Just to simplify the example, we consider the following formula $\Phi:=x=2s\otimes(y_1,y_2)\wedge y_1\setminus y_2=y_1$ therefore the
      $2s\otimes$-nodes that we need to build the model for $\Phi$ are $Q^{(2s\otimes,\Sigma,\mathfrak{I})}(y_1,y_2)$.

      Then $Q^{(2s\otimes,\Sigma,\mathfrak{I})}(y_1,y_2)=\{N_1=\{\sigma_1',\sigma_2\},
     N_2=\{\sigma_1'',\sigma_2\}\}$.

 $\Sigma_{2s\otimes}=\Sigma^{*}=\{\sigma_5',\sigma_5''\}$ is a $2s\otimes$-subpartition of $\Sigma$ since $\bigcup\Sigma^{*}$ is equal to $2s\otimes N_1\cup 2s\otimes N_2=2s\otimes (\mathfrak{I}(y_1),\mathfrak{I}(y_2))$.

 Since $\bigcup \Sigma_{2s\otimes}=2s\otimes N_1\cup 2s\otimes N_2 $,  $\sub{\Pi}{C}=\{N_1,N_2\}$ which are $2s\otimes$-blocks (in particular $2s\otimes$-upblocks). Then $\mathfrak{I}(x)=2s\otimes(\mathfrak{I}(y_1),\mathfrak{I}(y_2))$

The vertex set of $\sub{\mathcal{G}}{\Sigma}$ are
$\{\sigma_1',\sigma_1'',\sigma_2,\sigma_5',\sigma_5'',N_1,N_2\}$, membership edges of $\sub{\mathcal{G}}{\Sigma}$ are $\{\langle \sigma_1',N_1 \rangle,\langle \sigma_2,N_1 \rangle,\langle \sigma_1'',N_2 \rangle,\langle \sigma_2,N_2 \rangle\}$ and
saturation edges of $\sub{\mathcal{G}}{\Sigma}$ are $\{\langle N_1, \sigma_5'\rangle,\langle N_2,\sigma_5'' \rangle\}$.

In a rather natural way an $ns\otimes$ operator ({\bf N Singlelizer operator})  could be defined.
\end{myexample}

\subsubsection{Accessible $C$-graphs}

The following definitions of source places and accessible $C$-graphs are the obvious generalization to $C$ of those in \cite{CU21} Definition 5.
\begin{mydef}[Accessible $C$-graphs]\label{def:accessibility}\rm
A place of a $C$-graph $\mathcal{G}=(\Places,\Nodes,\TARGETS)$ is a \textsc{source place} if it has no incoming edges. The remaining places, namely those with incoming edges, are called \textsc{$C$-places}. We denote by $\sub{\Places}{C}$ the collection of the $C$-places of $\mathcal{G}$.

A place of $\mathcal{G}$ is \textsc{accessible} (from the source places of $\mathcal{G}$) if either it is a source place or, recursively, it is the target of some node of $\mathcal{G}$ whose places are \emph{all} accessible from the source places of $\mathcal{G}$. Finally, a $C$-graph is \textsc{accessible} when all its places are accessible.\footnote{Thus, a $C$-graph with no source places is trivially not accessible.}
\end{mydef}
A $C$-graph $\sub{\mathcal{G}}{\Sigma}$ induced by a given partition $\Sigma$ is always accessible, the proof is straightforward and is left to the reader.

Only \emph{accessible} $C$-graphs are relevant for our decidability purposes.

Again following \cite{CU21} just interchanging $\otimes$ with $\mathcal{C}$,
the accessible $C$-graph induced by a partition $\Sigma$ satisfying a given $\BST$ $\mathcal{C}$-conjunction $\myphi$ \emph{fulfills} $\myphi$, according to the following definition.
{\color{red}nel seguito per la definizione di Q $\mathfrak{F}$ svolgera il ruolo di $\mathfrak{I}$ e P il ruolo di $Sigma$}
\begin{mydef}[Fulfillment by an accessible $C$-graph]\label{satAccessible}\rm
An accessible $C$-graph $\mathcal{G}=(\Places,\Nodes,\TARGETS)$  \textsc{fulfills} a given \BST$\mathcal{C}$-conjunction $\myphi$ provided that there exists a map $\mathfrak{F} \colon \Vars{\myphi} \rightarrow \pow{\Places}$ (called a \textsc{$\mathcal{G}$-fulfilling map for $\myphi$})  such that the following conditions are satisfied:
\begin{enumerate}[label=(\alph*), ref=\alph*]
\item\label{satAccessibleA} $\mathfrak{F}(x) = \mathfrak{F}(y) \star \mathfrak{F}(z)$, for every conjunct $x=y \star z$ in $\myphi$, where $\star \in \{\cup,\setminus\}$;

\item\label{satAccessibleB} $\mathfrak{F}(x) \neq \mathfrak{F}(y)$, for every conjunct $x \neq y$ in $\myphi$;

\item\label{satAccessibleC} for every conjunct $x=\mathcal{C}(x_1\dots x_n)$ in $\myphi$,

\begin{enumerate}[label=(\ref{satAccessibleC}$_{\arabic*}$)]
\item\label{c1} $\emptyset \neq \Targets{N}\subseteq \mathfrak{F}(x)$, for all $N\in Q_{x_1,\cdots,x_n}$;

\item\label{c2} $\mathfrak{F}(x) \subseteq \bigcup_{N\in Q_{x_1,\cdots,x_n}}\TARGETS[N]$;

\item\label{c3} $\bigcup\TARGETS[\Nodes \setminus Q_{x_1,\cdots,x_n}] \cap \mathfrak{F}(x) = \emptyset$.
\end{enumerate}
\end{enumerate}

\end{mydef}
\begin{myexample}
  The $C$-graph shown in the Example \ref{cgraph} fulfills $\myphi$ of the same example.
\end{myexample}

\begin{mylemma}\label{C-GraphFulfillsMyphi}
The accessible $C$-graph induced by a partition satisfying a given \BST$\mathcal{C}$-conjunction $\myphi$ fulfills $\myphi$.
\end{mylemma}

The proof of the above lemma follows that one of \cite{CU21} using the property of operators which constructors.
For reader's convenience we report it in the appendix.
As an immediate consequence, we have:

\begin{mycorollary}\label{cor-otimesGraphFulfillsMyphi}
A satisfiable $\BST\mathcal{C}$-conjunction with $n$ variables is fulfilled by an accessible $C$-graph of size (at most) $2^{n}-1$.
\end{mycorollary}

\section{Decidability of Extensions of $\BST$ with $H$-operators}

We investigate operators which constructors that allow filling any accessible $C$-graph. This goal is achieved through proper set assignment to sources.
As described above, places in $\Places$ are intended to be an abstract representation of the blocks of a partition $\Sigma$. We will fill these abstract blocks with sets. Whenever we do this, we point out the match between the abstract block and the block filled of sets via a bijection $q \mapsto q^{(\bullet)}$ from $\Places$ onto $\Sigma$. The image of this bijection is the partition we will create.

Assume $C$ is a set-constructor and $\Places$ defined as before. Let $H_0$ be the $H_0=\bigcup_{q\in S_0}q^{(\bullet)}$, $S_0\subseteq\Places$. Let $C(B_0)=\bigcup_{q\in S_1}q^{(\bullet)}$, where $B_0\subseteq S_0$ and $S_1\subseteq\Places\setminus S_0$.

We continue inductively by defining $C(B_{i})=\bigcup_{q\in S_{i+1}}q^{(\bullet)}$, where $B_i\subseteq S_0\cup\cdots\cup S_i$ and $S_{i+1}\subseteq\Places\setminus (S_0\cup\cdots\cup S_i)$.

{\color{red}Da sistemare}
\begin{mydef}\label{Hoper}

A set-constructor $C$ is an $H$-constructor if for any set of places $\Places$ there exists a finite set $H_0$ of any cardinality such that

\begin{enumerate}[label=(\alph*), ref=\alph*]
\item\label{HConstructor1} for all $i$ $C(B_i)\cap H_0=\emptyset$;

\item\label{HConstructor2} for any $q\in B$, $\card{C(B)}\ge \card{q}$.

\end{enumerate}

An operator $\mathcal{C}$ with constructors $(C,Q)$ is an $H$-operator if $C$ is an $H$-constructor.

\end{mydef}

\begin{myremark}
  Examples of $H$-operators are, for instance, those described in Example \ref{constructors}.
  Indeed, in all three cases, $H_0$ can be composed of sets of the same rank. Trivially, there exists a collection of sets of this kind of any cardinality.
  Since they are them all disjoint operator and they increase the rank of the sets with each application of them, they meet the requirements to be $H$-operators.
\end{myremark}

\subsection{Construction process}
The Lemma~\ref{C-GraphFulfillsMyphi} can be reversed.

\begin{mylemma}\label{rev-otimesGraphFulfillsMyphi}
If a \BST$\mathcal{C}$-conjunction is fulfilled by an accessible $C$-graph, then it is satisfiable.
\end{mylemma}

Observe that the combination of Lemma \ref{C-GraphFulfillsMyphi} and Lemma \ref{rev-otimesGraphFulfillsMyphi} implies the following Theorem
\begin{mytheorem}\label{Decidability}
  Any language of the type \BST$C$ where $C$ is an H-operator is decidable.
\end{mytheorem}

The proof of Lemma \ref{rev-otimesGraphFulfillsMyphi} runs exactly as the analogue in \cite{CU21} except for the use of the properties of $H$-operator $\mathcal{C}$. Indeed, in the initialization we use the set $H_0$ to fill source places.
Then, when you have to fill the remaining places on the graph, you use the constructor $C$ on the nodes $Q$. To be sure that you build a real partition you have to be sure that all the elements created by the constructor $C$ have not been created before.
These elements cannot
belong to $H_0$ by property (a) of the $H$-operator, and cannot belong to
other places, $C$ being a disjoint operator (by Definition \ref{disjoint}). This ensures that what you are building
it is a partition.
The fact that $H_0$ can be of any cardinality together with property (\ref{HConstructor2}) of $H$-operator ensure that you have enough elements to follow all paths of the graph.

The remaining details of the proof are the same.
For the sake of completeness, we report the whole proof in the Appendix.
Moreover, just substituting $\otimes$ with $\mathcal{C}$, the proof of finite satisfiability runs in the same way.

\section{An exponential obstruction to polynomial-size explicit $C$-graphs}
\label{sec:graph-obstruction}

The graph-based satisfiability criterion of the preceding sections suggests
searching for a polynomial-size explicit accessible $C$-graph that fulfills
every satisfiable $\BST\mathcal C$-conjunction. Such a small-graph
property would give an $NP$ upper bound whenever fulfillment and
accessibility of the proposed graph can be checked in polynomial time.
We prove here that the required small-graph property is false \emph{already
for} $\mathcal C=\otimes$. Thus the specific pruning strategy that keeps
only polynomially many places and redirects constructor edges cannot
supply a general $NP$-membership proof.

Throughout this section,
\[
  A\otimes B:=\{\{a,b\}:a\in A,\ b\in B\}.
\]
The singleton $\{a\}$ is included as the case $a=b$. We use only
literals $x=y\cup z$, $x=y\setminus z$, $x=y\otimes z$ and $x\ne y$,
which are admissible in $\BST^\otimes$.

\subsection{A family of quadratic-size formulas}

Fix $m\ge1$. Introduce variables
\[
 U,Z,\quad E_1,\ldots,E_m,\quad
 X_1,\ldots,X_m,\quad Y_1,\ldots,Y_m.
\]
Let $\Phi_m$ be the conjunction of the following literals:
\begin{align}
 Z&=U\setminus U,                                      \label{eq:bad-zero}\\
 E_i&\ne Z,                                            &&1\le i\le m,\label{eq:bad-nonempty}\\
 Z&=E_i\setminus U,\qquad Z=X_i\setminus U,          &&1\le i\le m,\label{eq:bad-subset}\\
 Y_i&=X_i\cup E_i,\qquad X_i=Y_i\otimes U,            &&1\le i\le m,\label{eq:bad-product}\\
 E_i&=E_i\setminus E_j,                               &&1\le i\ne j\le m,\label{eq:bad-disjointE}\\
 E_i&=E_i\setminus X_j,                               &&1\le i,j\le m.\label{eq:bad-disjointX}
\end{align}
In particular, $Z=\varnothing$, $E_i\ne\varnothing$,
$E_i\cup X_i\subseteq U$, the $E_i$ are pairwise disjoint,
and $E_i\cap X_j=\varnothing$ for every $i,j$.
There are $3m+2$ variables and $O(m^2)$ literals; in a conventional
binary encoding of indexed variable names the length is
$O(m^2\log(m+1))$.

\begin{mylemma}\label{lem:bad-satisfiable}
For every $m\ge1$, the conjunction $\Phi_m$ is satisfiable in the
universe of well-founded sets.
\end{mylemma}
\begin{proof}
Let $U=\mathrm{HF}$ be the set of all hereditarily finite sets.
Choose pairwise distinct $e_1,\ldots,e_m\in\mathrm{HF}$ such that
$|e_i|\ge3$ for all $i$, and put $E_i=\{e_i\}$ and $Z=\varnothing$.
For every $i$ define an increasing sequence by
\[
 A_i^{(0)}=\varnothing,\qquad
 A_i^{(r+1)}=(A_i^{(r)}\cup\{e_i\})\otimes U
 \quad(r\ge0),
\]
and set
\[
 X_i=\bigcup_{r<\omega}A_i^{(r)},\qquad
 Y_i=X_i\cup E_i.
\]
Because the unordered pair of any two hereditarily finite sets is
hereditarily finite, $A_i^{(r)}\subseteq U$ for every $r$ and
$X_i\subseteq U$. The operation $B\mapsto B\otimes U$ preserves
increasing unions, so
\[
 X_i=(X_i\cup E_i)\otimes U=Y_i\otimes U.
\]
Each member of any $X_i$ has cardinality one or two, whereas
$|e_j|\ge3$; hence $e_j\notin X_i$ for all $i,j$. The sets $E_i$
are pairwise disjoint, nonempty and contained in $U$. All literals
\eqref{eq:bad-zero}--\eqref{eq:bad-disjointX} follow.
\end{proof}

\subsection{Exponential growth of Venn profiles}

For a set assignment $M$ satisfying $\Phi_m$ and an element
$a\in M(U)$, let
\[
 \operatorname{prof}_M(a)
 :=\{i\in[m]:a\in M(X_i)\},\qquad [m]=\{1,\ldots,m\}.
\]
\begin{mylemma}\label{lem:bad-venn}
Every model $M\models\Phi_m$ realizes every nonempty profile
$S\subseteq[m]$: there exists $w_S\in M(U)$ with
\[
  \operatorname{prof}_M(w_S)=S.
\]
Consequently its Venn partition with respect to the $X_i$ has at least
$2^m-1$ nonempty regions.
\end{mylemma}
\begin{proof}
Choose $e_i\in M(E_i)$ for each $i$, using
\eqref{eq:bad-nonempty}. By \eqref{eq:bad-subset}--\eqref{eq:bad-disjointX},
$e_i\in M(U)$ and
\[
 e_i\in M(Y_j)\quad\Longleftrightarrow\quad i=j.
 \tag{$*$}
\]
For each singleton $\{i\}$, put $w_{\{i\}}=\{e_i\}$. Since
$e_i\in M(Y_i)\cap M(U)$,
$w_{\{i\}}\in M(X_i)=M(Y_i)\otimes M(U)$.
On the other hand, membership of $\{e_i\}$ in
$M(Y_j)\otimes M(U)$ forces $e_i\in M(Y_j)$; by $(*)$ this is
impossible when $j\ne i$. Thus
$\operatorname{prof}_M(w_{\{i\}})=\{i\}$, and
$w_{\{i\}}\in M(U)$ by \eqref{eq:bad-subset}.

Assume that $w_S$ has been constructed for a nonempty proper
$S\subset[m]$, and choose $i\notin S$. Define
\[
 w_{S\cup\{i\}}=\{w_S,e_i\}.
\]
The two entries belong to $M(U)$. For every $j\in S$,
$w_S\in M(X_j)\subseteq M(Y_j)$, so
$w_{S\cup\{i\}}\in M(Y_j)\otimes M(U)=M(X_j)$.
The same holds for $j=i$ because $e_i\in M(Y_i)$.
If $j\notin S\cup\{i\}$, then
$e_i\notin M(Y_j)$ by $(*)$, and
$w_S\notin M(Y_j)$: the inductive profile excludes $M(X_j)$,
whereas $w_S$ belongs to some $M(X_k)$ with $k\in S$,
and \eqref{eq:bad-disjointX} excludes membership in $M(E_j)$.
Therefore $w_{S\cup\{i\}}\notin M(Y_j)\otimes M(U)=M(X_j)$.
We have obtained the required profile, and
$w_{S\cup\{i\}}\in M(U)$. Induction completes the proof.
Distinct profiles determine distinct nonempty Venn regions.
\end{proof}

\subsection{The lower bound for every fulfilling $\otimes$-graph}

The preceding lemma concerns actual models. For the complexity argument
we need the stronger assertion that the growth is forced by the
\emph{abstract fulfillment conditions} themselves. In the terminology of
Definition~\ref{satAccessible}, for $x=y\otimes z$ write
\[
 Q_{y,z}=\big\{\{p,r\}:p\in\mathfrak F(y),\ r\in\mathfrak F(z)\big\}.
\]
We use conditions (c$_1$) and (c$_3$) of that definition:
for every node $N\in Q_{y,z}$ its target set is nonempty and lies in
$\mathfrak F(x)$, whereas if $N\notin Q_{y,z}$ none of its
targets lies in $\mathfrak F(x)$.

\begin{mytheorem}[Exponential graph lower bound]
\label{thm:bad-graphs}
For every $m\ge1$, $\Phi_m$ is satisfiable and each standard
fulfilling accessible $\otimes$-graph
$\mathcal G=(\Places,\Nodes,\TARGETS)$ for $\Phi_m$ satisfies
\[
 \boxed{\ |\Places|\ge 2^m-1\ }.
\]
In fact the lower bound only uses the fulfillment conditions,
not accessibility.
\end{mytheorem}
\begin{proof}
Satisfiability follows from Lemma~\ref{lem:bad-satisfiable}.
Let $\mathfrak F$ be a fulfilling map for an arbitrary such graph and put
\[
 \chi(q):=\{j\in[m]:q\in\mathfrak F(X_j)\}
 \qquad(q\in\Places).
\]
The Boolean fulfillment conditions yield
\begin{align*}
 \mathfrak F(Z)&=\varnothing,
 &\mathfrak F(E_i)&\ne\varnothing,\\
 \mathfrak F(E_i)\cup\mathfrak F(X_i)&\subseteq\mathfrak F(U),
 &\mathfrak F(Y_i)&=\mathfrak F(X_i)\cup\mathfrak F(E_i),
\end{align*}
and
\[
 \mathfrak F(E_i)\cap\mathfrak F(E_j)=\varnothing\ (i\ne j),
 \qquad
 \mathfrak F(E_i)\cap\mathfrak F(X_j)=\varnothing\quad(\forall i,j).
\]
Choose a place $p_i\in\mathfrak F(E_i)$. It belongs to
$\mathfrak F(U)$ and $\mathfrak F(Y_i)$, but to no
$\mathfrak F(Y_j)$ with $j\ne i$.
For each $j$ set
\[
 Q_j:=Q_{Y_j,U}
 =\big\{\{p,r\}:p\in\mathfrak F(Y_j),\
 r\in\mathfrak F(U)\big\}.
\]
Then $N_i:=\{p_i\}\in Q_i$ and $N_i\notin Q_j$ for $j\ne i$.
Condition (c$_1$) applied to the $i$th product literal gives
$\TARGETS(N_i)\ne\varnothing$. Choose
$q_{\{i\}}\in\TARGETS(N_i)$. By (c$_1$) for $i$ and
(c$_3$) for all $j\ne i$ we obtain
\[
 \chi(q_{\{i\}})=\{i\}.
\]
Furthermore $q_{\{i\}}\in\mathfrak F(U)$, and it belongs to no
$\mathfrak F(E_j)$.

Assume that for a nonempty $S\subset[m]$ there exists
$q_S\in\mathfrak F(U)$, belonging to no $\mathfrak F(E_j)$, with
$\chi(q_S)=S$. Fix $i\notin S$ and put
$N:=\{q_S,p_i\}$. Both places belong to $\mathfrak F(U)$, and
\[
 N\in Q_j\quad\Longleftrightarrow\quad j\in S\cup\{i\}.
\]
Since $N\in Q_i$, its target set is nonempty. Choose
$q_{S\cup\{i\}}\in\TARGETS(N)$. Applying (c$_1$) to all
$j\in S\cup\{i\}$ and (c$_3$) to all other $j$, we have
\[
 \chi(q_{S\cup\{i\}})=S\cup\{i\}.
\]
As the new place belongs to $\mathfrak F(X_i)$, the Boolean
conditions also imply that it lies in $\mathfrak F(U)$ and in no
$\mathfrak F(E_j)$. Hence induction produces a place $q_S$ for
every nonempty $S\subseteq[m]$. Such places are pairwise distinct
because their $\chi$-profiles are distinct. Consequently
$|\Places|\ge2^m-1$.
\end{proof}

\begin{myproposition}[An explicit accessible graph for $\Phi_m$]
\label{prop:bad-explicit-graph}
For each $m\ge1$ there exists a fulfilling accessible $\otimes$-graph
for $\Phi_m$ with $m+2^m-1$ places.
\end{myproposition}
\begin{proof}
Let
\[
 \Places=\{b_i:i\in[m]\}\ \dot\cup\
 \{q_S:\varnothing\ne S\subseteq[m]\},
\]
and define $\lambda(b_i)=\{i\}$ and $\lambda(q_S)=S$.
Take as $\otimes$-nodes all singleton and doubleton subsets
$N=\{p,r\}\subseteq\Places$, allowing $p=r$; the other
subsets of $\Places$ have no outgoing arrows. Put
\[
 \TARGETS(\{p,r\})=\{q_{\lambda(p)\cup\lambda(r)}\}.
\]
Every such target is well defined since the labels are nonempty.
Set
\[
 \begin{aligned}
 \mathfrak F(U)&=\Places,&\mathfrak F(Z)&=\varnothing,\\
 \mathfrak F(E_i)&=\{b_i\},&
 \mathfrak F(X_i)&=\{q_S:i\in S\},\\
 &&\mathfrak F(Y_i)&=\mathfrak F(X_i)\cup\{b_i\}.
 \end{aligned}
\]
All Boolean conditions of $\Phi_m$ hold. Since both entries of
any product node lie in $\mathfrak F(U)$, for each $i$
\[
 \{p,r\}\in Q_{Y_i,U}
 \quad\Longleftrightarrow\quad
 i\in\lambda(p)\cup\lambda(r)
 \quad\Longleftrightarrow\quad
 \TARGETS(\{p,r\})\subseteq\mathfrak F(X_i).
\]
This gives (c$_1$) and (c$_3$); moreover, for $q_S\in\mathfrak F(X_i)$
with $i\in S$, the node $\{b_i,q_S\}$ belongs to $Q_{Y_i,U}$
and has target $q_S$, proving (c$_2$).
Each $b_i$ is a source. The place $q_{\{i\}}$ is reached from
$\{b_i\}$, and for $|S|\ge2$, $q_S$ is reached from
$\{b_i,q_{S\setminus\{i\}}\}$ with $i\in S$.
Therefore every place is accessible, and the graph fulfills $\Phi_m$.
\end{proof}

\begin{mycorollary}[Failure of polynomial-size explicit graph witnesses]
\label{cor:no-poly-graphs}
There exists no polynomial $p$ such that every satisfiable
$\BST^\otimes$-conjunction $\Phi$ has a standard fulfilling accessible
$\otimes$-graph with at most $p(|\Phi|)$ places. In particular, the
polynomial-place pruning step of the previous $NP$-membership argument
cannot hold for all $\BST^\otimes$-conjunctions, and therefore cannot
prove such an upper bound uniformly for the class of $H$-operators
containing $\otimes$.
\end{mycorollary}
\begin{proof}
The formulas $\Phi_m$ have length $O(m^2\log(m+1))$, while
Theorem~\ref{thm:bad-graphs} requires at least $2^m-1$ places for
every fulfilling graph. This exceeds any polynomial in the formula
length for sufficiently large $m$.
\end{proof}

\begin{myremark}[Precise scope of the obstruction]
The corollary rules out a \emph{particular type of certificate}:
an explicitly listed fulfilling accessible $\otimes$-graph with
polynomially many places. It neither proves
$\BST^\otimes\notin NP$ nor rules out succinct encodings of large
fulfilling graphs with polynomial-time verification, or other kinds
of certificates. Since ordinary $\BST$ is a fragment of
$\BST^\otimes$, the satisfiability problem for $\BST^\otimes$
is $NP$-hard; the matching $NP$ upper bound is not obtained by the
argument refuted here. For a restricted fragment, Cantone and Maugeri
\cite{CM23Tensor} prove $NP$-completeness when the number of distinct
left-hand-side variables in product literals is logarithmic in the
formula length.
\end{myremark}

\begin{myremark}[A separate issue with injectivity]
The unordered Cartesian product is not injective on all set arguments
in the sense previously assumed in the small-graph proof: for any
$A\ne B$, one has $A\otimes\varnothing=B\otimes\varnothing
=\varnothing$. Thus injectivity without restrictions on empty
arguments cannot justify redirecting graph edges either. The
exponential lower bound above is independent of this observation.
\end{myremark}

\section*{Appendix}

\textbf{Proof of Lemma \ref{C-GraphFulfillsMyphi}}
\begin{proof}

Let us assume that $\Sigma$ satisfies a conjunction $\myphi$ via a partition assignment $\mathfrak{I} \colon \Vars{\myphi} \rightarrow \pow{\Sigma}$.

The proof of properties (\ref{satAccessibleA}),(\ref{satAccessibleB}) of Definition \ref{satAccessible} are plain verifications due to the fact that the assignment $\mathfrak{I}$ satisfies $\myphi$.

Regarding (\ref{satAccessibleC}) of Definition \ref{satAccessible},
let $x = \mathcal{C}(x_1\cdots x_n)$ be a $\mathcal{C}$-literal in $\myphi$, so that
\begin{equation}\label{C-identity}
\medcup \mathfrak{I}(x) = \mathcal{C}(\medcup \mathfrak{I}(x_1),\dots ,\medcup\mathfrak{I}(x_n))
\end{equation}
holds.

We preliminarily observe that $\mathfrak{I}(x) \subseteq \sub{\Sigma}{C}$.

Indeed, by \eqref{C-identity}, $\medcup \mathfrak{I}(x) = \mathcal{C}(\medcup \mathfrak{I}(x_1),\dots ,\medcup\mathfrak{I}(x_n))$ since $\mathcal{C}$ is an operator which constructors
$\bigcup_{N\in Q(x_1\cdots x_n)} C(N)=\mathcal{C}(\medcup\mathfrak{I}(x_1)\cdots\medcup\mathfrak{I}(x_n))$, and therefore $\mathfrak{I}(x)$ is a $C$-subpartition of $\Sigma$.

Let $\upsilon_1 \in \sub{\mathfrak{F}}{\Sigma}(x_1)\dots\upsilon_n \in \sub{\mathfrak{F}}{\Sigma}(x_n)$ ($\mathfrak{F}$ is the syntactical analogue of $\mathfrak{I}$, that is
$
\mathfrak{F}(x) \defAs \{ q \in \Places^* \st q^{(\bullet)} \in \mathfrak{I}(x)\}\/, \quad \text{for $x \in \Vars{\myphi}$}
$
). Then $\upsilon_1^{(\bullet)} \subseteq \medcup \mathfrak{I}(x_1)$ and $\upsilon_n^{(\bullet)} \subseteq \medcup \mathfrak{I}(x_n)$, and consequently, by \eqref{C-identity}, $C(\upsilon_1^{(\bullet)}, \dots ,\upsilon_n^{(\bullet)}) \subseteq \medcup \mathfrak{I}(x)$.
Hence,
\begin{align*}
\emptyset &\neq \medcup \big\{ q^{(\bullet)} \in \Sigma \st q^{(\bullet)} \cap C(\upsilon_1^{(\bullet)}, \dots ,\upsilon_n^{(\bullet)}) \neq \emptyset \big\} \subseteq \medcup \mathfrak{I}(x)\\
\intertext{and therefore}
\emptyset &\neq \big\{ q^{(\bullet)} \in \Sigma \st q^{(\bullet)} \cap C(\upsilon_1^{(\bullet)}, \dots ,\upsilon_n^{(\bullet)}) \neq \emptyset \big\} \subseteq  \mathfrak{I}(x).\\
\intertext{Thus,}
\emptyset &\neq  \big\{ q \in \sub{\Places}{\Sigma} \st q \cap C(\upsilon_1, \dots ,\upsilon_n) \neq \emptyset \big\} \subseteq \sub{\mathfrak{F}}{\Sigma}(x) \cap \sub{\Places}{\Sigma ,C},\\
\intertext{so that}
\emptyset &\neq \sub{\TARGETS}{\Sigma}(\{\upsilon_1,\dots ,\upsilon_n \}) = \big\{ q \in \sub{\Places}{\Sigma,\!C} \st q \cap C(\upsilon_1, \dots ,\upsilon_n) \neq \emptyset \big\} \subseteq \sub{\mathfrak{F}}{\Sigma}(x).
\end{align*}

From arbitrariness of $q,\upsilon_1,\dots ,\upsilon_n$, we get property \ref{c1} of Definition \ref{satAccessible}.

Let $q \in \sub{\mathfrak{F}}{\Sigma}(x)$, so that $q^{(\bullet)} \subseteq \medcup \mathfrak{I}(x)$.
Let $s \in q^{(\bullet)}$. Hence, $s \in C(\upsilon_1^{(\bullet)}, \dots ,\upsilon_n^{(\bullet)})$, for some
$\upsilon_1,\dots ,\upsilon_n \in \sub{\Places}{\Sigma}$ such that $\upsilon_1^{(\bullet)} \in \mathfrak{I}(x_1),\dots ,\upsilon_n^{(\bullet)} \in \mathfrak{I}(x_n)$.

Since $q^{(\bullet)} \cap C(\upsilon_1^{(\bullet)}, \dots ,\upsilon_n^{(\bullet)}) \neq \emptyset$ and $q^{(\bullet)} \in \sub{\Sigma}{C}$, we have
$
q \in \sub{\TARGETS}{\Sigma}(\{\upsilon_1,\dots ,\upsilon_n\}) \subseteq \medcup_{B\in Q_{x_1,\dots,x_n}}\sub{\TARGETS}{\Sigma}(B).
$
Hence, from the arbitrariness of $q \in \sub{\mathfrak{F}}{\Sigma}(x)$, it follows that
\begin{equation}\label{FSigmax-inclusion}
\sub{\mathfrak{F}}{\Sigma}(x) \subseteq \medcup_{B\in Q_{x_1,\dots,x_n}} \sub{\TARGETS}{\Sigma}(B).
\end{equation}
which proves property \ref{c2} of Definition \ref{satAccessible}.

Finally we have to prove property \ref{c3} of Definition \ref{satAccessible}:
\begin{equation}\label{emptyIntersection}
\medcup \sub{\TARGETS}{\Sigma}[\sub{\Nodes}{\Sigma}] \setminus \medcup_{B\in Q_{x_1,\dots,x_n}}\sub{\TARGETS}{\Sigma}(B) \cap \sub{\mathfrak{F}}{\Sigma}(x) = \emptyset.
\end{equation}
By contradiction, let us assume that there exists some $q \in \medcup \sub{\TARGETS}{\Sigma}[\sub{\Nodes}{\Sigma}] \setminus \medcup_{B\in Q_{x_1,\dots,x_n}}\sub{\TARGETS}{\Sigma}(B) \cap \sub{\mathfrak{F}}{\Sigma}(x)$.
Hence, $q \in \sub{\TARGETS}{\Sigma}(A)$ for some $A \in \sub{\Nodes}{\Sigma}$, therefore
$q^{(\bullet)} \cap C(A^{(\bullet)}) \neq \emptyset$. By \eqref{C-identity},  $q^{(\bullet)} \subseteq \medcup \mathfrak{I}(x)$, then $q^{(\bullet)}\in \mathfrak{I}(x)$, therefore $q\in\mathfrak{F}(x)$ which is a contradiction by \ref{c2}, thus proving \eqref{emptyIntersection}.

We can summarize what we have just proved by saying that the accessible $C$-graph induced by a partition $\Sigma$ satisfying a given \BST$\mathcal{C}$-conjunction $\myphi$ \emph{fulfills} $\myphi$.
\end{proof}

\medskip

\textbf{Proof of Lemma \ref{rev-otimesGraphFulfillsMyphi}}

\begin{proof}
Let $\mathcal{G}=(\Places,\Nodes,\TARGETS)$ be an accessible $C$-graph, and let us assume that $\mathcal{G}$ fulfills a given $\BST$ $\mathcal{C}$-conjunction $\myphi$ via the map $\mathfrak{F} \colon \Vars{\myphi} \rightarrow \pow{\Places}$.

To each place $q \in \Places$, we associate a set $q^{(\bullet)}$, initially empty. Then, by suitably exploiting the $C$-graph $\mathcal{G}$ as a kind of flow graph, we shall show that the sets $q^{(\bullet)}$'s can be monotonically extended by a (possibly infinite) \emph{construction process} (comprising a finite \emph{initialization phase} and a subsequent (possibly infinite) \emph{stabilization phase}) in such a way that the following properties hold:
\begin{enumerate}[label=(P\arabic*)]
\item\label{P1} After each step, the sets $q^{(\bullet)}$'s are pairwise disjoint.

\item\label{P2}
At the end of the initialization phase all the $q^{(\bullet)}$'s are nonempty (and pairwise disjoint). Thus, after each step in the subsequent stabilization phase,  the sets $q^{(\bullet)}$'s, with  $q \in \Places$, form a partition equipollent with $\Places$.

\item\label{P3}
After each step in the stabilization phase, the inclusion
\[
q^{(\bullet)} \subseteq \medcup \big\{ C(A^{(\bullet)}) \st A \in \TARGETS^{-1}(q)\big\}
\]
holds, for each $C$-place $q \in \sub{\Places}{C}$, where we are using the notation $B^{(\bullet)} \defAs \{p^{(\bullet)} \st p \in B\}$ for $B \in \Nodes$.
\end{enumerate}

\begin{enumerate}[resume*]
\item\label{P4} At the end of the construction process, we have
\[
C(A^{(\bullet)}) \subseteq \bigcup \{q^{(\bullet)} \st q \in \Targets{A}\},
\]
for each $C$-node $A \in \sub{\Nodes}{C}$ (namely for each node $A \in \Nodes$ such that $\Targets{A} \neq \emptyset$).\footnote{Should the construction process involve denumerably many steps, the final values of the $q^{(\bullet)}$'s are to be intended as limit of the sequences of their values after each step in the stabilization phase.}
\end{enumerate}

Subsequently, we shall prove that the properties \ref{P1}--\ref{P4} together with the conditions \eqref{satAccessibleA}--\eqref{satAccessibleC} of Definition~\ref{satAccessible}, characterizing the fulfilling $C$-graph $\mathcal{G}$, allow one to show that the partition $\{q^{(\bullet)} \st q \in \Places\}$ resulting from the above construction process satisfies our conjunction $\myphi$.

\bigskip

The initialization and stabilization phases of our construction process consist of the following steps.
\paragraph{Initialization phase:}
\begin{enumerate}[label=(I$_{\arabic*}$)]
\item\label{I1} To begin with, let $\{\overline q \st q \in \Places \setminus \sub{\Places}{C}\}$ be any partition equipollent to the set $\Places \setminus \sub{\Places}{C}$ of the source places of $\mathcal{G}$, where each block $\overline q$, for $q \in \Places \setminus \sub{\Places}{C}$, is filled by elements of the set $H_0$, that we choose of cardinality greater than $\card{\Places}^{k+1}$ where $k$ is the length of the longest path of the $C$-graph. This allow us to distribute at least $\card{\Places}^{k}$ to each source place, then put
\[
q^{(\bullet)} \defAs \begin{cases}
\overline q & \text{if } q \in \Places \setminus \sub{\Places}{C}\\
\emptyset & \text{if } q \in \sub{\Places}{C}.
\end{cases}
\]
\end{enumerate}

We say that a place $q \in \Places$ has already been \emph{initialized} when $q^{(\bullet)} \neq \emptyset$. Likewise, a $C$-node $A \in \sub{\Nodes}{C}$ has been \emph{initialized} when its places have all been initialized. During the initialization phase, an initialized $C$-node $A \in \sub{\Nodes}{C}$ is said to be \textsc{ready} if it has some target that has not been yet initialized.

\begin{enumerate}[resume*]
\item\label{I2} While there are places in $\Places$ not yet initialized, pick any ready node $A \in \Nodes$ and distribute evenly all the members of $C(A^{(\bullet)})$ among all of its targets.
\end{enumerate}

The accessibility of $\mathcal{G}$ guarantees that the while-loop \ref{I2} terminates in a finite number of iterations.

Concerning property \ref{P1}, we observe that since $C$ is a disjoint constructor so whenever it applies at each distribution step to different nodes it produces disjoint sets that are to distributed to the respective targets, moreover $\mathcal{C}$ is an $H$ operator therefore its constructor $C$ guarantees that $C(A^{(\bullet)})\cap H_0=\emptyset$.

At the end of the initialization phase all the $q^{(\bullet)}$'s are nonempty, so  property \ref{P2} holds. Indeed, if there were no $C$-places, then all places would be initialized just after step \ref{I1}, and so all the $q^{(\bullet)}$'s would be nonempty.
On the other hand, if $|\sub{\Places}{C}|>0$
consider any path $L$ of the $C$-graph.
The first arrow of $L$ start from a node $A$ composed of $\Places \setminus \sub{\Places}{C}$ of the source places of $\mathcal{G}$. Each of the places has at least $\card{\Places}^{k}$ of elements.
Since $\mathcal{C}$ is an $H$-operator $C(A^{(\bullet)})\ge\card{q^{(\bullet)}}\ge\card{\Places}^{k}$.
Hence, each of the $|\Targets{A}| \leq |\sub{\Places}{C}|$ sets $t^{(\bullet)}$, for $t \in \Targets{A}$, will receive at least $\card{\Places}^{k-1}$ elements by the distribution step relative to the node $A$.

Inductively, at the $h$-arrow of the path $L$ each target has at least $\card{\Places}^{k-h}$ elements and this guarantees that all path can be followed until the end, then at the end of the while-loop \ref{I2} we shall have $|q^{(\bullet)}|\neq\emptyset$, for each $q \in \Places$.

\paragraph{Stabilization phase:} During the stabilization phase, a $C$-node $A \in \sub{\Nodes}{C}$ is \textsc{ripe} if
\[
C(A^{(\bullet)}) \setminus \medcup \big\{q^{(\bullet)} \st q \in \Targets{A}\big\} \neq \emptyset.
\]

We execute the following (possibly infinite) loop:
\begin{enumerate}[label=(S$_{\arabic*}$)]
\item\label{while-loop} While there are ripe $C$-nodes, pick any of them, say $A \in \Nodes$, and distribute all the members of $C(A^{(\bullet)}) \setminus \medcup \big\{q^{(\bullet)} \st q \in \Targets{A}\big\}$ (namely the members of $C(A^{(\bullet)})$ that have not been distributed yet) among its targets howsoever.
\end{enumerate}

The fairness condition that one must comply with is the following:
\begin{quote}
once a $C$-node becomes ripe during the stabilization phase, it must be picked for distribution within a finite number of iterations of the while-loop \ref{while-loop}.
\end{quote}
A possible way to enforce such condition consists, for instance, in maintaining all ripe $C$-nodes in a queue $\mathcal{Q}$, picking always the $C$-node to be used in a distribution step from the front of $\mathcal{Q}$ and adding the $C$-nodes that have just become ripe to the back of $\mathcal{Q}$, provided that they are not already in $\mathcal{Q}$.

By induction on $n \in \Nats$, it is not hard to show that properties \ref{P1} and \ref{P3} will hold just after the $n$-th iteration of the while-loop \ref{while-loop} of the stabilization phase, and that property \ref{P4} will hold at the end of the stabilization phase, in case of termination.

Instead, when the stabilization phase runs for denumerably many steps, the final partition $\Places^{(\bullet)}$ is to be intended as the limit of the partial partitions constructed after each step of the stabilization phase. Specifically, for each place $q \in \Places$, we let $q^{(i)}$ be the value of $q^{(\bullet)}$ just after the $i$-th iteration of \ref{while-loop}. Plainly, we have
\begin{equation}\label{monotonicity}
q^{(i)} \subseteq q^{(i+1)}, \quad \text{for } i \in \Nats.
\end{equation}
Then we put
\begin{equation}\label{q-bullet-limit}
q^{(\bullet)} \defAs \bigcup_{i \in \Nats} q^{(i)}, \quad \text{for } q \in \Places
\end{equation}
(notation overloading should not be a problem).

By way of illustration, we prove that property \ref{P4} holds for the partition $\Places^{(\bullet)} = \big\{q^{(\bullet)} \st q \in \Places \big\}$, when the $q^{(\bullet)}$'s are defined by \eqref{q-bullet-limit}. To this purpose, let $A \in \Places$ be such that $\Targets{A} \neq \emptyset$, and assume for contradiction that
\[
C(A^{(\bullet)}) \not\subseteq \medcup \big\{q^{(\bullet)} \st q \in \Targets{A}\big\}.
\]
Let $s$ be any element in $C(A^{(\bullet)}) \setminus \medcup \big\{q^{(\bullet)} \st q \in \Targets{A}\big\}$, and let $i \in \Nats$ be the smallest index such that $s \in C(A^{(i)})$, where $A^{(i)} \defAs \{q^{(i)} \st q \in A \}$. Since $s \in C(A^{(i)}) \setminus \medcup \big\{q^{(i)} \st q \in \Targets{A}\big\}$, the node $A$ must have been ripe just after the $i$-th iteration of \ref{while-loop}. Therefore, by the fairness condition, the node $A$ will be picked for distribution in a finite number of steps, say $k$, after the $i$-th step, so that we have
\begin{align*}
C(A^{(i)}) &\subseteq C(A^{(i+k)}) && \text{(by \eqref{monotonicity})}\\
&\subseteq \medcup \big\{q^{(i+k+1)} \st q \in \Targets{A}\big\}\\
&\subseteq \medcup \big\{q^{(\bullet)} \st q \in \Targets{A}\big\},
\end{align*}
and therefore $s \in \medcup \big\{q^{(\bullet)} \st q \in \Targets{A}\big\}$, which is a contradiction. Thus, property \ref{P4} holds also when the construction process takes a denumerable number of steps.

\bigskip

Next, we show that the final partition $\Places^{(\bullet)} = \{q^{(\bullet)} \st q \in \Places \}$ satisfies $\myphi$. In particular, we prove that the partition assignment $\mathfrak{I} \colon \Vars{\myphi} \rightarrow \pow{\Places^{(\bullet)}}$ defined by
\[
\mathfrak{I}(x) \defAs \{q^{(\bullet)} \st q \in \mathfrak{F}(x) \} \text{, \quad for } x \in \Vars{\myphi},
\]
satisfies $\myphi$, where we recall that $\mathfrak{F}$ is the $\mathcal{G}$-fulfilling map for $\myphi$.

Since $\mathfrak{F}$ is a $\mathcal{G}$-fulfilling map for $\myphi$, then
\begin{enumerate}[label=-]
\item for every literal $x=y \star z$ in $\myphi$, with $\star \in \{\cup,\setminus\}$, we have $\mathfrak{F}(x) = \mathfrak{F}(y) \star \mathfrak{F}(z)$, so that $\mathfrak{I}(x) = \mathfrak{I}(y) \star \mathfrak{I}(z)$ holds; and

\item for every literal $x \neq y$ in $\myphi$, we have $\mathfrak{F}(x) \neq \mathfrak{F}(y)$, so that $\mathfrak{I}(x) \neq \mathfrak{I}(y)$ holds.
\end{enumerate}
Thus, by Lemma~\ref{partitionAssignmentBoolean}, the partition assignment $\mathfrak{I}$ satisfies all Boolean literals in $\myphi$ of types
\[
x=y \cup z, \quad x=y \setminus z, \quad x \neq y.
\]

Next, let $x=C(x_1\cdots x_n)$ be a conjunct of $\myphi$. We prove separately that the following inclusions hold:
\begin{gather}
\medcup \mathfrak{I}(x) \subseteq \mathcal{C}(\medcup \mathfrak{I}(x_1),\cdots ,\medcup \mathfrak{I}(x_n)) \label{firstInclusion}\\
\mathcal{C}(\medcup \mathfrak{I}(x_1),\cdots ,\medcup \mathfrak{I}(x_n)) \subseteq \medcup \mathfrak{I}(x).\label{secondInclusion}
\end{gather}

We recall that $\mathcal{C}$ is an operator with constructor then

$\bigcup_{N\in Q(x_1\cdots x_n)} C(N^{(\bullet)})=\mathcal{C}(\medcup\mathfrak{I}(x_1)\cdots\medcup\mathfrak{I}(x_n))$

Concerning \eqref{firstInclusion}, let $q^{(\bullet)} \subseteq \bigcup\mathfrak{I}(x)$. Then $q^{(\bullet)} \in \mathfrak{I}(x)$, so that $q \in \mathfrak{F}$. By \ref{c2}, $q$ cannot be a source place. Hence, by \ref{P3}, we have:
\[
q^{(\bullet)} \subseteq \medcup \big\{ C(A^{(\bullet)}) \st A \in \TARGETS^{-1}(q)\big\}.
\]

By definition \ref{satAccessible} (\ref{satAccessibleB}) and (\ref{satAccessibleC}) $\TARGETS^{-1}(q)\subseteq Q(x_1\cdots x_n)$ therefore
\[
q^{(\bullet)} \subseteq \bigcup_{N\in Q(x_1\cdots x_n)}C(N^{(\bullet)}).
\]
which in turn implies
\[
\bigcup_{q\in\mathfrak{I}(x)}q^{(\bullet)} \subseteq \mathcal{C}(\medcup\mathfrak{I}(x_1)\cdots\medcup\mathfrak{I}(x_n)).
\]
and finally implies \ref{firstInclusion}.

Concerning the inclusion \eqref{secondInclusion}, let $s \in \mathcal{C}(\medcup \mathfrak{I}(x_1),\cdots, \medcup \mathfrak{I}(x_n))$. Hence, $s \in C(N^{(\bullet)})$, for some $N \in Q_{x_1,\cdots,x_n}$.

From \ref{c1}, we have $\emptyset \neq \Targets{N} \subseteq \mathfrak{F}(x)$. Thus, by \ref{P4},
\begin{align*}
C(N^{(\bullet)}) &\subseteq \medcup \{ q^{(\bullet)} \st q \in \Targets{N} \}\\
&\subseteq \medcup \{ q^{(\bullet)} \st q \in  \mathfrak{F}(x) \}\\
&= \medcup \mathfrak{I}(x),
\end{align*}
and therefore $s \in \medcup \mathfrak{I}(x)$, proving \eqref{secondInclusion} by the arbitrariness of $s \in \mathcal{C}(\medcup \mathfrak{I}(x_1),\cdots, \medcup \mathfrak{I}(x_n))$.

Hence, the partition assignment $\mathfrak{I}$ satisfies also all the literals in $\myphi$ of the form $x=\mathcal{C}(x_1,\cdots,x_n)$, and in turn the final partition $\Places^{(\bullet)}$ satisfies the conjunction $\myphi$.
\end{proof}

By combining Lemmas~\ref{C-GraphFulfillsMyphi} and \ref{rev-otimesGraphFulfillsMyphi} and Corollary~\ref{cor-otimesGraphFulfillsMyphi}, we obtain:
\begin{mytheorem}
A \BST$\mathcal{C}$-conjunction with $n$ variables is satisfiable if and only if it is fulfilled by an accessible $C$-graph of size (at most) $2^{n}-1$.
\end{mytheorem}


\end{document}